\documentclass[10pt]{amsart}
\usepackage{graphicx}
\usepackage{xcolor}
\usepackage{amssymb}
\usepackage{amsmath}
\usepackage[T1]{fontenc}
\newtheorem{theorem}{Theorem}[section]
\newtheorem{lemma}[theorem]{Lemma}

\newtheorem{remark}[theorem]{Remark}
\newtheorem{definition}[theorem]{Definition}
\newenvironment{proof of main thm 1}{{\bf Proof of Theorem \ref{main thm 1}.}}{\hfill\fbox{}\par\vspace{.2cm}}
\newenvironment{proof of main thm 2}{{\bf Proof of Theorem \ref{main thm 2}.}}{\hfill\fbox{}\par\vspace{.2cm}}
\newenvironment{proof of main thm 3}{{\bf Proof of Theorem \ref{main thm 3}.}}{\hfill\fbox{}\par\vspace{.2cm}}
\newenvironment{proof of main thm 4}{{\bf Proof of Theorem \ref{main thm 4}.}}{\hfill\fbox{}\par\vspace{.2cm}}
\newenvironment{proof of thm}{{\bf Proof.}}{\hfill\fbox{}\par\vspace{.2cm}}
\newenvironment{proof of main three}{{\bf Proof of Theorem \ref{main three}.}}{\hfill\fbox{}\par\vspace{.2cm}}
\newenvironment{proof of sub thm 2}{{\bf Proof of Theorem \ref{sub thm 2}.}}{\hfill\fbox{}\par\vspace{.2cm}}
\newenvironment{proof of sub thm 3}{{\bf Proof of Theorem \ref{sub thm 3}.}}{\hfill\fbox{}\par\vspace{.2cm}}
\numberwithin{equation}{section}

\def\charf {\mbox{{\text 1}\kern-.24em {\text l}}}

\def\bea{\begin{eqnarray*}}
\def\eea{\end{eqnarray*}}
\def\be{\begin{eqnarray}}
\def\ee{\end{eqnarray}}

\begin{document}

\title[Results for $\mathfrak{g}$-stable with boundary]{Some results on the $\mathfrak{g}$-stability of surfaces with boundary}



\author{Sanghun Lee}
\address{Department of Mathematics and Institute of Mathematical Science, Pusan National University, Busan 46241, Korea}
\email{kazauye@pusan.ac.kr}


\thanks{The author was supported by the National Research Foundation of Korea (NRF) grant funded by the Korea government (No. RS-2023-00214167).}


\date{}

\dedicatory{}

\begin{abstract}
In this paper, we investigate the geometric properties associated with the $\mathfrak{g}$-stability of surfaces with boundary whose null expansion satisfies $\Theta^{+} = h \geq 0$. First, we show that a $\mathfrak{g}$-stable hypersurface with free boundary admits a metric of positive scalar curvature with minimal boundary under suitable conditions. Second, for $\mathfrak{g}$-stable surfaces with free boundary, we derive an area estimate and determine the topology of the surface. Finally, we extend our free boundary results to the case of capillary boundary.
\end{abstract}

\maketitle

\section{Introduction}

An \textit{initial data set} $(M^{n+1}, g, A^{M})$ consists of an $(n+1)$-dimensional manifold $M$ that arise as a spacelike hypersurface in an $(n+2)$-dimensional spacetime $(\mathsf{L}^{n+2}, \bar{g})$, where $A^{M}$ denotes the second fundamental form of $M$. 

In initial data sets, minimal surfaces in Riemannian geometry are replaced by \textit{marginally outer trapped surfaces} (\textit{MOTS}s). The notion of a MOTS was introduced in the development of black hole theory, and MOTSs serve as geometric objects that link spacetime geometry with the physics of black holes. For example, under suitable conditions, MOTSs arise as cross sections of event horizons in stationary black hole spacetimes. 

In pure mathematics, MOTSs appeared in the work of Schoen and Yau \cite{SY4} in connection with Jang’s equation and the positive mass theorem. Since then, many results on minimal surfaces in Riemannian geometry have been extended to the spacetime using MOTSs (see \cite{ALY, ALM, AMS, AM2, CG2, CH, EHLS, GA1, GM2, LLU}). 

Recently, Galloway and Mendes \cite{GM1} extended basic results concerning closed MOTSs to the general null expansion. A detailed description of general null expansion is provided in Section 2.

In this paper, we further extend results on closed surfaces with general null expansion \cite{GM1, LBS} to surfaces with boundary.

Before stating our results, we introduce the following notions. Let the \textit{momentum tensor} $\pi$ be defined by $\pi = A^{M} - \left({\rm tr}_{g}A^{M}\right)g$. We denote by $(\iota_{\bar{N}}\pi)^{T}$ the restriction of $\iota_{\bar{N}}\pi = \pi(\bar{N},\cdot)$ to tangent vector fields of $\partial M$ in $M$, where $\bar{N}$ is the outward unit normal of $\partial M$ in $M$. We say that $(M, g, A^{M})$ satisfies the \textit{dominant boundary energy condition} (\textit{DBEC}) if 
\begin{align} \label{DBEC}
H^{\partial M} \geq \vert \left(\iota_{\bar{N}}\pi\right)^{T} \vert
\end{align}
(see \cite{ALM, AM, M2}).

Our first result extends Theorem 3.1 of \cite{GM1} and Theorem 1 of \cite{LBS} to the setting of hypersurfaces with free boundary. This result is also related to the existence of a metric of positive scalar curvature with minimal boundary on an $n$-dimensional hypersurface $\Sigma$.

\begin{theorem} \label{main thm 1}
Let $(M^{n+1}, g, A^{M})$ be an $(n+1)$-dimensional initial data set with boundary, and let $\Sigma^{n}$ be a compact $n$-dimensional $\mathfrak{g}$-stable hypersurface with free boundary such that the null expansion satisfies $\Theta^{+} = h \in C^{\infty}(\Sigma)$, with $h \geq 0$. Suppose that $(M^{n+1}, g, A^{M})$ satisfies $\mu - \vert J \vert \geq C_{0}$ for some constant $C_{0} \in \mathbb{R}^{+} \cup \{0\}$ and the DBEC. If one of the following conditions holds:
\begin{itemize}
\item[(i)] $C_{0} = 0$ and $\tau \leq \frac{h}{2}$ along $\Sigma$, $h \not\equiv 0$, or
\item[(ii)] $C_{0} > 0$ and $h\tau \leq C_{0}$ along $\Sigma$,
\end{itemize}
then $\Sigma$ admits a metric of positive scalar curvature with minimal boundary.
\end{theorem}

For three-dimensional initial data sets with boundary, the rigidity results for stable MOTS ($\Theta^{+} = 0$) with free boundary were established in \cite{ALY} (see also \cite{AM, M2}). In this paper, we extend these results to the setting of $\mathfrak{g}$-stability.

\begin{theorem} \label{main thm 2}
Let $(M^{3}, g, A^{M})$ be a three-dimensional initial data set with boundary satisfying $\mu - \vert J \vert \geq C_{1}$ for some positive constant $C_{1}$ and the DBEC. Suppose that $\Sigma^{2}$ is a compact $\mathfrak{g}$-stable surface with free boundary such that $\Theta^{+} = h \in C^{\infty}(\Sigma)$, $h \geq 0$, and $\tau \leq \frac{h}{2}$. Then
\begin{align} \label{area free}
A(\Sigma) \leq \frac{2\pi}{C_{1}}
\end{align}
and $\Sigma$ is topologically a two-disk. Furthermore, if equality holds in (\ref{area free}), then the following properties holds:
\begin{itemize}
    \item[(1)] The null second fundamental form $\chi^{+}$ vanishes on $\Sigma$. In particular, $\Sigma$ is a MOTS with free boundary.
    \item[(2)] $\mu + J(N) = \mu - \vert J \vert = C_{1}$ on $\Sigma$, and $H^{\partial M} = \vert (\iota_{\bar{N}}\pi)^{T} \vert$ along $\partial\Sigma$.
    \item[(3)] The first eigenvalue $\lambda_{1}(L)$ is equal to zero.
    \item[(4)] $\Sigma$ has constant Gaussian curvature $K^{\Sigma} = C_{1}$, and $\partial\Sigma$ has zero geodesic curvature $k^{\partial\Sigma} = 0$.
\end{itemize}
\end{theorem}

Next, we introduce a generalized notion of the DBEC for the case where $\Sigma$ has a capillary boundary. This condition is called the \textit{tilted dominant boundary energy condition} (\textit{TDBEC}) is defined by 
\begin{align} \label{TDBEC}
H^{\partial M} + (\cos\theta){\rm tr}_{\partial M}A^{M} \geq \sin\theta\vert A^{M}(\bar{N}, \cdot)^{T} \vert, 
\end{align}
where $A^{M}(\bar{N},\cdot)^{T}$ denotes the tangential component of the 1-form $A^{M}(\bar{N},\cdot)$ along $\partial M$ (see \cite{CH, CW}). In the free boundary case ($\theta = \frac{\pi}{2}$), this condition reduces to the DBEC.

We now extend our result for the free boundary (Theorem \ref{main thm 1}) to the setting of capillary boundary.

\begin{theorem} \label{main thm 3}
Let $(M^{n+1}, g , A^{M})$ be an $(n+1)$-dimensional initial data set with boundary, and let $\Sigma^{n}$ be a compact $n$-dimensional $\mathfrak{g}$-stable hypersurface with capillary boundary such that $\Theta^{+} = h \in C^{\infty}(\Sigma)$, $h \geq 0$, and $\theta \in [\frac{\pi}{2}, \pi)$. Suppose that $(M^{n+1}, g, A^{M})$ satisfies $\mu - \vert J \vert \geq C_{0}$ for some constant $C_{0} \in \mathbb{R}^{+} \cup \{0\}$ and the TDBEC. If either
\begin{itemize}
\item[(i)] $C_{0} = 0$ and $\tau \leq \frac{h}{2}$ along $\Sigma$, $h \not\equiv 0$, or
\item[(ii)] $C_{0} > 0$ and $h\tau \leq C_{0}$ along $\Sigma$,
\end{itemize}
then $\Sigma$ admits a positive scalar curvature with minimal boundary.
\end{theorem}
 
Finally, we establish a rigidity result for $\mathfrak{g}$-stable surfaces with capillary boundary in three-dimensional initial data sets with boundary.

\begin{theorem} \label{main thm 4}
Let $(M^{3}, g, A^{M})$ be a three-dimensional initial data set with boundary satisfying $\mu - \vert J \vert \geq C_{1}$ for some positive constant $C_{1}$ and the TDBEC. Suppose that $\Sigma^{2}$ is a compact $\mathfrak{g}$-stable surface with capillary boundary such that $\Theta^{+} = h \in C^{\infty}(\Sigma)$, $h \geq 0$, $\tau \leq \frac{h}{2}$, and $\theta \in [\frac{\pi}{2}, \pi)$. Then
\begin{align} \label{area capillary}
A(\Sigma) \leq \frac{2\pi}{C_{1}}
\end{align}
and $\Sigma$ is topologically a two-disk. Moreover, if equality holds in (\ref{area capillary}), then the following properties hold:
\begin{itemize}
    \item[(1)] The null second fundamental form satisfies $\chi^{+} = 0$ on $\Sigma$; in particular, $\Sigma$ is a MOTS with capillary boundary.
    \item[(2)] $\mu + J(N) = \mu - \vert J \vert = C_{1}$ on $\Sigma$, and $H^{\partial M} + (\cos\theta){\rm tr}_{\partial M}A^{M} = \sin\theta \vert A^{M}(\bar{N}, \cdot)^{T} \vert$ along $\partial\Sigma$.
    \item[(3)] The first eigenvalue satisfies $\lambda_{1}(L_{1}) = 0$ on $\Sigma$.
    \item[(4)] $\Sigma$ has constant Gaussian curvature $K^{\Sigma} = C_{1}$, and $\partial\Sigma$ has zero geodesic curvature $k^{\partial\Sigma} = 0$.
\end{itemize}
\end{theorem}

\section{Preliminaries}

Let $(\mathsf{L}^{n+2}, \bar{g})$ be an $(n+2)$-dimensional spacetime satisfying the Einstein equation
\begin{align*} \label{EM}
G_{\bar{g}} \equiv Ric_{\bar{g}} - \frac{1}{2}R_{\bar{g}}\bar{g} = T,
\end{align*}
where $Ric_{\bar{g}}$ and $R_{\bar{g}}$ denote the Ricci curvature and scalar curvature of $\bar{g}$, respectively, and $T$ is a symmetric 2-tensor.

An $(n+1)$-dimensional initial data set $(M^{n+1}, g, A^{M})$ is a spacelike hypersurface in $(\mathsf{L}^{n+2}, \bar{g})$ endowed with the induced Riemannian metric $g$ and the second fundamental form $A^{M}$ of $M$ in $\mathsf{L}$. More precisely, the second fundamental form $A^{M}$ is defined by 
\begin{align*}
A^{M}(X,Y) = \bar{g}(\nabla^{\mathsf{L}}_{X}N_{M}, Y),
\end{align*}
where $N_{M}$ is the future-directed timelike unit normal vector to $M$, and $X,Y \in T_{p}M$.

By the Gauss-Codazzi equations, an initial data set $(M^{n+1}, g, A^{M})$ satisfies the Einstein constraint equations
\begin{align*}
&\mu := T(N_{M},N_{M}) = \frac{1}{2}\left(R^{M} + \tau^{2} - \vert A^{M} \vert^{2}\right), \\
&J := T(N_{M}, \cdot) = {\rm div}_{g}\left(A^{M} - \tau g\right),
 \end{align*}
where $R^{M}$ denotes the scalar curvature of $M$, and $\tau = {\rm tr}_{g}A^{M}$ is the mean curvature of $M$ in $\mathsf{L}$. We refer to $\mu$ as the local \textit{energy density} and to $J$ as the local \textit{current density}. The \textit{dominant energy conditions} then requires that
\begin{align*}
\mu \geq \vert J \vert \,\, {\rm  on} \,\, M.
\end{align*}

Let $\Sigma^{n} \subset M^{n+1}$ be an embedded $n$-dimensional hypersurface, and let $N$ denote a unit normal vector field to $\Sigma$ in $M$. By convention, we choose $N$ to be the outward pointing unit normal, so that $-N$ is inward pointing. Then $\Sigma$ admits two future-directed null normal vector fields given by
 \begin{align*}
\ell^{\pm} = N_{M} \pm N.
 \end{align*}
We define the \textit{null second fundamental form} $\chi^{\pm}$ by
\begin{align*}
\chi^{\pm}(X,Y) = \bar{g}(\nabla^{\mathsf{L}}_{X}\ell^{\pm}, Y) = A^{M}|_{\Sigma} \pm A^{\Sigma},
\end{align*}
where $X, Y \in T\Sigma$, and $A^{\Sigma}$ denotes the second fundamental form of $\Sigma$ in $M$. The \textit{null expansion scalars} $\Theta^{\pm}$ of $\Sigma$ in $M$ with respect to $N$ are defined by
\begin{align*}
\Theta^{\pm} = {\rm div}_{\Sigma}\ell^{\pm} = {\rm tr}_{\Sigma}A^{M} \pm H,
\end{align*}
where $H$ is the mean curvature of $\Sigma$ in $M$. When $\Theta^{+} = 0$, we say that $\Sigma$ is a \textit{marginally outer trapped surface} (\textit{MOTS}). In the time-symmetric case ($A^{M} = 0$), a MOTS reduces to a minimal surface in Riemannian geometry.
 
We now introduce the notion of \textit{$\mathfrak{g}$-stability} using \textit{general null expansion} $\Theta^{+} = h$, which generalizes the usual notion of stability for MOTS. This notion was introduced by Galloway and Mendes in \cite{GM1} and is slightly different from the stability for prescribed null expansion considered in \cite{CH, CH2}. For sufficiently small $\epsilon_{0} > 0$, the map
\begin{align*}
\Psi:[0, \epsilon_{0}] \times \Sigma \rightarrow M,
\end{align*}
defined by $\Psi(t,p) = {\rm exp}_{p}(tN(p))$, is well-defined. On the other hand, given a positive function $\varphi \in C^{\infty}(\Sigma)$ with $\Vert \varphi \Vert_{C^{0}} \leq \epsilon_{0}$, the map
\begin{align*}
\Psi^{\varphi}: [0,1] \times \Sigma \rightarrow M,
\end{align*}
defined by $\Psi^{\varphi}(t,p) = {\rm exp}_{p}(t\varphi(p)N(p))$, is also well-defined. We denote $\Psi^{\varphi}(t,\Sigma)$ by $\Sigma^{\varphi}_{t}$ and we write $\Theta^{\pm}_{\varphi}(t)$ for the null expansion scalars of $\Sigma^{\varphi}_{t}$.

\begin{definition} [Definition 2.1 in \cite{GM1}]
 We say that $\Sigma$ is \textit{$\mathfrak{g}$}-stable if
\begin{align*}
\frac{\partial \Theta^{+}_{\varphi}}{\partial t}\Big{|}_{t = 0} \geq 0
\end{align*}
for some function $\varphi > 0$ with $\Vert \varphi \Vert_{C^{0}} \leq \epsilon_{0}$.
\end{definition}
It is well known that (see \cite{AMS, AM})
\begin{align*}
\frac{\partial\Theta^{+}_{\varphi}}{\partial t}\Big{|}_{t = 0} = L\varphi,
\end{align*}
where $L: C^{\infty}(\Sigma) \rightarrow C^{\infty}(\Sigma)$ is the elliptic operator defined by
\begin{align}
L\varphi = -\Delta_{\Sigma}\varphi + 2\langle X, \nabla^{\Sigma}\varphi \rangle + \left(Q - \vert X \vert^{2} + {\rm div}_{\Sigma}X - \frac{h^{2}}{2} + h\tau\right)\varphi,
\end{align}
where $h = \Theta^{+}_{\varphi}(0)$ is the null expansion of $\Sigma^{\varphi}_{0} = \Sigma$, $\tau = {\rm tr}_{g}A^{M}$ is the mean curvature of $M$, and
\begin{align*}
Q = \frac{R^{\Sigma}}{2} - \left(\mu + J(N)\right) - \frac{\vert \chi^{+} \vert^{2}}{2}.
\end{align*}
Here $R^{\Sigma}$ denotes the scalar curvature of $\Sigma$, and $X$ is the vector field tangent to $\Sigma$ that is dual to the $1$-form $A^{M}(N,\cdot)|_{\Sigma}$.

It is well-known that there exists a real number $\lambda$, called the \textit{principal eigenvalue} of $L$, such that $L\varphi = \lambda\varphi$ for some positive eigenfunction $\varphi \in C^{\infty}(\Sigma)$, and that $\lambda \leq {\rm Re}(\eta)$ for any other eigenvalue $\eta$ of $L$ (see \cite{AMS}). The eigenspace of $L$ associated with $\lambda$ is one-dimensional. In particular, $\Sigma$ is $\mathfrak{g}$-stable if and only if $\lambda \geq 0$ (see \cite{GM1}). We denote the principal eigenvalue $\lambda$ of $L$ by $\lambda_{1}(L)$.

We now introduce the notion of \textit{$\mathfrak{g}$-stability with capillary boundary} for initial data set with boundary. We begin by recalling the definition of a \textit{hypersurface with capillary boundary}. Let $(M, g, A^{M})$ be an $(n+1)$-dimensional initial data set with boundary, and let $\Sigma$ be a compact \textit{properly} embedded hypersurface with boundary in $M$, meaning that $\Sigma$ is embedded in $M$ and $\partial\Sigma = \Sigma \cap \partial M$. We say that $\Sigma$ has a \textit{capillary boundary} in $M$ if $\Sigma$ meets $\partial M$ at a constant contact angle $\theta \in (0,\pi)$. More precisely, assume that $\Sigma$ separates $M$ into two components, and label one of them by $\Omega$. Let $N$ be the unit normal vector to $\Sigma$ pointing into $\Omega$. Let $\bar{\nu}$ denote the unit normal to $\partial\Sigma$, pointing out of $\Omega$, as a subset of $\partial M$. Furthermore, let $\nu$ be the outward unit normal to $\partial\Sigma$ in $\Sigma$, and let $\bar{N}$ be the outward unit normal of $\partial M$. Then one can choose orthonormal bases $\{N, \nu \}$ and $\{\bar{N}, \bar{\nu} \}$ such that the following relations hold (cf. \cite{RS}):
\begin{align*}
&\bar{N} = (\cos\theta) N + (\sin\theta)\nu, \\
&\bar{\nu} = -(\sin\theta)N + (\cos\theta)\nu.
\end{align*}
Equivalently,
\begin{align*}
&N = -(\sin\theta)\bar{\nu} + (\cos\theta)\bar{N}, \\
&\nu = (\cos\theta)\bar{\nu} + (\sin\theta)\bar{N},
\end{align*}
where $\theta$ is the constant contact angle between $\Sigma$ and $\partial M$.

Next, we introduce the notion of a \textit{$\mathfrak{g}$-stable hypersurface with capillary boundary}. Let $(\Sigma_{t})_{\vert t \vert < \epsilon}$ be a smooth variation of $\Sigma$, and consider the functional
\begin{align*}
F[\Sigma_{t}] = \int_{\Sigma_{t}} \left(\Theta^{+}(t) - h \right) \langle W_{t}, N_{t} \rangle dv + \int_{\partial\Sigma_{t}}\langle W_{t}, \nu_{t} - ({\rm cos} \, \theta)\bar{\nu}_{t} \rangle ds,
\end{align*}
where the variational vector field is given by $W_{t} = \frac{\partial}{\partial t}$ and $h \in C^{\infty}(\Sigma)$ (see also \cite{ALY, M2}). Since the general null expansion satisfies $\Theta^{+} = h$, the first variation of $F[\Sigma_{t}]$ is given by
\begin{align} \label{eigen ftnc}
\frac{\partial}{\partial t}\Bigr|_{t = 0}F[\Sigma_{t}] =& \int_{\Sigma}-\varphi\Delta_{\Sigma}\varphi + 2\varphi\langle X, \nabla^{\Sigma}\varphi \rangle + \left(Q + {\rm div}_{\Sigma}X - \vert X \vert^{2} + h\tau - \frac{h^{2}}{2}\right)\varphi^{2} dv \\
&+ \int_{\partial\Sigma}\left(\varphi\frac{\partial\varphi}{\partial \nu} - q\varphi^{2}\right), \nonumber
\end{align}
where $\varphi \in C^{\infty}(\Sigma)$ and $q = -(\cot\theta) A^{\Sigma}(\nu,\nu) + \frac{1}{\sin\theta}A^{\partial M}(\bar{\nu},\bar{\nu})$. The eigenvalue problem associated with the functional (\ref{eigen ftnc}) is then given by
\begin{align} \label{eigen capillary}
\begin{cases}
&L\varphi := -\Delta_{\Sigma}\varphi + 2\langle X,\nabla^{\Sigma}\varphi \rangle + \left(Q + {\rm div}_{\Sigma}X - \vert X \vert^{2} + h\tau - \frac{h^{2}}{2} \right)\varphi = \lambda \varphi \quad {\rm on} \,\, \Sigma, \\
&B_{c}\varphi := \frac{\partial\varphi}{\partial\nu} - q\varphi = 0 \quad {\rm along} \,\, \partial\Sigma.
\end{cases}
\end{align}
\begin{definition} \label{stable with capillary}
We say that $\Sigma$ is a \textit{$\mathfrak{g}$-stable hypersurface with capillary boundary} if there exists a non-negative function $\varphi \in C^{\infty}(\Sigma)$, $\varphi \not\equiv 0$, satisfying Robin boundary condition $B_{c}\varphi = 0$ such that $L\varphi \geq 0$ on $\Sigma$.
\end{definition}
When the contact angle $\theta$ equals $\frac{\pi}{2}$, we say that $\Sigma$ has a \textit{free boundary}. More precisely, along $\partial\Sigma$, the unit normal vector $\nu$ of $\partial\Sigma$ in $\Sigma$ coincides with the normal vector $\bar{N}$ of $\partial M$. In the free boundary case, the associated eigenvalue problem is given by
\begin{align} \label{eigen free}
\begin{cases}
&L\varphi := -\Delta_{\Sigma}\varphi + 2\langle X,\nabla^{\Sigma}\varphi \rangle + \left(Q + {\rm div}_{\Sigma}X - \vert X \vert^{2} + h\tau - \frac{h^{2}}{2} \right)\varphi = \lambda \varphi \quad {\rm on} \,\, \Sigma, \\
&B_{f}\varphi := \frac{\partial\varphi}{\partial\nu} - A^{\partial M}(N,N)\varphi = 0 \quad {\rm along} \,\, \partial\Sigma.
\end{cases}
\end{align}
\begin{definition}
We say that $\Sigma$ is a $\mathfrak{g}$-stable hypersurface with free boundary if there exists a non-negative function $\varphi \in C^{\infty}(\Sigma)$, $\varphi \not\equiv 0$, satisfying the boundary condition $B_{f}\varphi = 0$ such that $L\varphi \geq 0$ on $\Sigma$.
\end{definition}

\begin{remark}
The notion of $\mathfrak{g}$-stability for a hypersurface extends, in the variational sense, the usual notion of stability for MOTSs with $\Theta^{+} = 0$.
\end{remark}

\section{Results for free boundary}

In this section, we establish several results concerning $\mathfrak{g}$-stable (hyper)surfaces with free boundary. First, we investigate conditions under which an $n$-dimensional $\mathfrak{g}$-stable hypersurface with free boundary admits a metric of positive scalar curvature and has minimal boundary. Before presenting the proof, we recall an important result that will be used in establishing our main theorem.

\begin{lemma} [Lemma 3.1 and Lemma 3.2 in \cite{M2}] \label{M2 thm}
Let $(M, g, A^{M})$ be an $(n+1)$-dimensional initial data set with boundary, and $\Sigma$ be a compact hypersurface with free boundary such that $\Theta^{+} = h \in C^{\infty}(\Sigma)$. Suppose that $\Sigma$ is a $\mathfrak{g}$-stable hypersurface with free boundary. Then the following hold:
\begin{itemize}
    \item[(1)] For all $\psi \in C^{\infty}(\Sigma)$ and $\tilde{Q} \in C^{\infty}(\Sigma)$, \begin{align} \label{thm 3.1}
0 \leq \int_{\Sigma} \psi\left(-\Delta_{\Sigma}\psi + \psi \tilde{Q}\right)dv + \int_{\partial\Sigma}\psi\left(\frac{\partial\psi}{\partial\nu} - \psi\left(A^{\partial M}(N,N) - \langle X, \nu \rangle\right)\right) ds.
    \end{align}
\end{itemize}
Hence, the first eigenvalue $\lambda_{1}(L_{1})$ of $L_{1} := -\Delta_{\Sigma} + \tilde{Q}$ on $\Sigma$, with Robin boundary condition $B_{1}\psi := B_{f}\psi + \langle X, \nu \rangle\psi = 0$ is non-negative. Moreover, if equality holds in (\ref{thm 3.1}) and $\psi \not\equiv 0 $, then $\psi$ is an eigenfunction of $L_{1}$ with $\lambda_{1}(L_{1}) = 0$.
\begin{itemize}
    \item[(2)] Let $\psi \in C^{\infty}(\Sigma)$ be positive and suppose that \begin{align*}
\begin{cases}
&\mathcal{L}\psi := -\Delta_{\Sigma}\psi + \left(\frac{R^{\Sigma}}{2} - \tilde{P}\right)\psi \geq 0 \quad {\rm on} \,\, \Sigma, \\
&\frac{\partial\psi}{\partial\nu} + H^{\partial\Sigma}\psi \geq 0 \quad {\rm along} \,\, \partial\Sigma,
\end{cases}
\end{align*}
\end{itemize}
where $R^{\Sigma}$ is the scalar curvature of $\Sigma$, $H^{\partial\Sigma}$ is the mean curvature of $\partial\Sigma$ in $\Sigma$, and $\tilde{P}$ is a non-negative function on $\Sigma$. Then $\Sigma$ admits a metric of positive scalar curvature with minimal boundary, unless $\Sigma$ is Ricci flat with totally geodesic boundary, $\tilde{P} = 0$, and $\psi$ is constant.
\end{lemma}

A detailed proof can be found in Section 3 of \cite{M2}.

We now prove our result on positive scalar curvature with minimal boundary.

\begin{proof of main thm 1}
Let $\psi > 0$ be an eigenfunction of $L_{1}\psi = -\Delta_{\Sigma}\psi + \tilde{Q}\psi$ on $\Sigma$ with Robin boundary condition $B_{1}\psi = 0$, associated with the first eigenvalue $\lambda_{1}(L_{1})$. Since $\Sigma$ is a $\mathfrak{g}$-stable hypersurface with free boundary, we obtain
\begin{align} \label{eq 3.2}
\begin{cases}
&L_{1}\psi = -\Delta_{\Sigma}\psi + \tilde{Q}\psi \geq 0 \,\, {\rm on} \,\, \Sigma, \\
&B_{1}\psi = \frac{\partial\psi}{\partial\nu} - \left(A^{\partial M}(N,N) - \langle X, \nu \rangle\right)\psi = 0 \,\, {\rm along} \,\, \partial\Sigma,
\end{cases}
\end{align}
Note that
\begin{align*}
\left(\iota_{\bar{N}}\pi\right)^{T}(N) = A^{M}(\bar{N}, N) - \tau g(\bar{N}, N) = \langle X, \bar{N} \rangle = \langle X, \nu \rangle \,\, {\rm along} \,\, \partial\Sigma.
\end{align*}
Therefore, we have
\begin{align} \label{eq 3.3}
A^{\partial M} (N, N) - \langle X, \nu \rangle &= H^{\partial M} - H^{\partial\Sigma} - \left(\iota_{\bar{N}}\pi\right)^{T}(N) \\
&\geq -H^{\partial\Sigma}, \nonumber
\end{align}
where, in the last inequality, we have used the DBEC \eqref{DBEC}. Substituting (\ref{eq 3.3}) into (\ref{eq 3.2}), we obtain
\begin{align*}
\begin{cases}
&-\Delta_{\Sigma}\psi + \tilde{Q}\psi \geq 0 \,\, {\rm on} \,\, \Sigma, \\
&\frac{\partial\psi}{\partial\nu} + H^{\partial\Sigma}\psi \geq 0 \,\, {\rm along} \,\, \partial\Sigma,
\end{cases}
\end{align*}
where $\tilde{Q} = \frac{R^{\Sigma}}{2} - \tilde{P}$ and $\tilde{P} = \left(\mu + J(N) \right) + \frac{\vert \chi^{+} \vert^{2}}{2} - h\left(\tau - \frac{h}{2}\right)$.

\begin{itemize}
    \item[CASE (i).] $C_{0} = 0$.
\end{itemize}
In this case, since $\mu + J(N) \geq \mu - \vert J \vert \geq 0$, we obtain
\begin{align*}
\tilde{P} &= \left(\mu + J(N)\right) + \frac{\vert \chi^{+} \vert^{2}}{2} - h\left(\tau - \frac{h}{2}\right) \\
&\geq -h\left(\tau - \frac{h}{2}\right) \geq 0,
\end{align*}
where, in the last inequality, we used the assumption $\tau \leq \frac{h}{2}$. If $\tilde{P} = 0$, then $0 = \mu + J(N) = h\left(\tau - \frac{h}{2}\right)$ and $\chi^{+} = 0$. This implies that $0 = {\rm tr}_{\Sigma}\chi^{+} = \Theta^{+} = h$, which contradicts the assumption that $h \not\equiv 0$. Therefore, $\Sigma$ admits a metric of positive scalar curvature with minimal boundary.

\begin{itemize}
    \item[CASE (ii).] $C_{0} > 0$.
\end{itemize}
In this case, we have
\begin{align*}
\tilde{P} &= \left(\mu + J(N)\right) + \frac{\vert \chi^{+} \vert^{2}}{2} - h\tau + \frac{h^{2}}{2} \\
&\geq C_{0} - h\tau \geq 0.
\end{align*}
If $\tilde{P} = 0$, then $C_{0} = \mu + J(N) = h\tau + \frac{h^{2}}{2}$. This implies that $\chi^{+} = 0$ and hence $h = 0$. On the other hand, since $C_{0} \geq h\tau = \mu + J(N) \geq C_{0}$, we must have $h\tau = C_{0} > 0$, which is a contradiction. Therefore, $\Sigma$ admits a metric of positive scalar curvature with minimal boundary.

\end{proof of main thm 1}

Next, we extend the rigidity results for stable MOTS ($\Theta^{+} = 0$) with free boundary to $\mathfrak{g}$-stable surface with free boundary under general null expansion ($\Theta^{+} = h$).

\begin{proof of main thm 2}
By (1) of Lemma \ref{M2 thm} and (\ref{eq 3.3}), we obtain
\begin{align*}
0 \leq \int_{\Sigma}\left(\vert \nabla^{\Sigma}\psi \vert^{2} + \tilde{Q}\psi^{2}\right) dv - \int_{\partial\Sigma}k^{\partial\Sigma}\psi^{2}ds,
\end{align*}
where $\tilde{Q} = \frac{R^{\Sigma}}{2} - \tilde{P}$ and $\tilde{P} = \left(\mu + J(N)\right) + \frac{\vert \chi^{+} \vert^{2}}{2} - h\left(\tau - \frac{h}{2}\right)$. Taking $\psi = 1$, we obtain
\begin{align*}
0 &\leq \int_{\Sigma} \tilde{Q} \, dv + \int_{\partial\Sigma} k^{\partial\Sigma} \, ds \\
&\leq \int_{\Sigma} K^{\Sigma} dv - \int_{\Sigma} \left(\mu - \vert J \vert\right)dv - \int_{\Sigma} \left(\tau - \frac{h}{2}\right)h \, dv + \int_{\partial\Sigma} k^{\partial\Sigma} \, ds \\
&\leq -C_{1}A(\Sigma) + 2\pi\chi(\Sigma),
\end{align*}
where $\chi(\Sigma)$ denotes the Euler characteristic of $\Sigma$. Hence, we conclude that
\begin{align} \label{eq 3.4.4}
C_{1}A(\Sigma) \leq 2\pi\chi(\Sigma) = 2\pi \left(2 - 2\alpha_{1} - \alpha_{2}\right),
\end{align}
where $\alpha_{1}$ is the genus of $\Sigma$ and $\alpha_{2}$ is the number of boundary components of $\Sigma$. This inequality implies that $\Sigma$ is topologically a disk. Moreover, if equality holds in (\ref{eq 3.4.4}), then all of the above inequalities become equalities. In particular:
\begin{itemize}
    \item[(1)] $\chi^{+} = 0$ and $\mu + J(N) = \mu - \vert J \vert = C_{1}$ on $\Sigma$;
    \item[(2)] $\psi = 1$ is an eigenfunction of $L_{1}$ associated with the first eigenvalue $\lambda_{1}(L_{1}) = 0$. Consequently,
    \begin{align*}
    \begin{cases}
    &0 = K^{\Sigma} - C_{1} \,\, {\rm on} \,\, \Sigma, \\
    &0 = - \left(A^{\partial M}(N,N) - \langle X,\nu \rangle\right) = k^{\partial\Sigma} \,\, {\rm along} \,\, \partial\Sigma;
    \end{cases}
    \end{align*}
    \item[(3)] hence, $\Sigma$ has constant Gaussian curvature equal to positive constant $C_{1}$ and its boundary has zero geodesic curvature.
\end{itemize}

\end{proof of main thm 2}

\section{Results for capillary boundary}

In this section, we extend our results for the free boundary case to the setting of capillary boundary. First, we present an analogous of Lemma \ref{M2 thm} for hypersurfaces with capillary boundary.

\begin{lemma} [Lemma 3.1 in \cite{LEE} and Lemma 3.2 in \cite{M2}] \label{thm 4.1}
Let $(M, g, A^{M})$ be an $(n+1)$-dimensional initial data set with boundary, and let $\Sigma$ be a compact hypersurface with capillary boundary such that $\Theta^{+} = h \in C^{\infty}(\Sigma)$. Suppose that $\Sigma$ is a $\mathfrak{g}$-stable hypersurface with capillary boundary. Then the following hold:
\begin{itemize}
    \item[(1)] For all $\psi \in C^{\infty}(\Sigma)$ and $\tilde{Q} \in C^{\infty}(\Sigma)$, \begin{align} \label{lem 4.1}
0 \leq \int_{\Sigma}\psi\left(-\Delta_{\Sigma}\psi + \psi\tilde{Q}\right)dv + \int_{\partial\Sigma}\psi\left(\frac{\partial\psi}{\partial\nu} - \psi \left(q - \langle X, \nu \rangle\right)\right) ds
    \end{align}
\end{itemize}
Hence, the first eigenvalue $\lambda_{1}(L_{1})$ of $L_{1} = -\Delta_{\Sigma} + \tilde{Q}$ on $\Sigma$ with Robin boundary condition $B_{2}\psi := B_{c}\psi + \langle X, \nu \rangle\psi = 0$ is non-negative. Moreover, if equality holds in (\ref{lem 4.1}) and $\psi \not\equiv 0$, then $\psi$ is an eigenfunction of $L_{1}$ associated with $\lambda_{1}(L_{1}) = 0$.
\begin{itemize}
    \item[(2)] Let $\psi \in C^{\infty}(\Sigma)$ be positive and suppose that \begin{align*}
\begin{cases}
&\mathcal{L}\psi := -\Delta_{\Sigma}\psi + \left(\frac{R^{\Sigma}}{2} - \tilde{P}\right)\psi \geq 0 \quad {\rm on} \,\, \Sigma, \\
&\frac{\partial\psi}{\partial\nu} + H^{\partial\Sigma}\psi \geq 0 \quad {\rm along} \,\, \partial\Sigma.
\end{cases}
\end{align*}
\end{itemize}
Then $\Sigma$ admits a metric of positive scalar curvature with minimal boundary, unless $\Sigma$ is Ricci flat with totally geodesic boundary, $\tilde{P} = 0$, and $\psi$ is constant.
\end{lemma}

The following lemma allows us to handle the boundary term $q$ more conveniently.

\begin{lemma}[Lemma 3.1 in \cite{LPP, LON}] \label{lem 4.2}
Let $(M, g, A^{M})$ be an $(n+1)$-dimensional initial data set with boundary $\partial M$, and let $\Sigma$ be a hypersurface with capillary boundary in $M$ with contact angle $\theta \in (0,\pi)$. Then
\begin{align*} 
A^{\partial M}(\bar{\nu},\bar{\nu}) - (\cos\theta)A^{\Sigma}(\nu,\nu) + (\sin\theta)H^{\partial\Sigma} = H^{\partial M} - (\cos\theta)H.  
\end{align*}
In particular, if $n = 2$, then $H^{\partial\Sigma} = k^{\partial\Sigma}$, where $k^{\partial\Sigma}$ denotes the geodesic curvature of $\partial\Sigma$ in $\Sigma$.
\end{lemma}

For a detailed proof, we refer the reader to Lemma 3.1 in \cite{LPP, LON}. By Lemma \ref{lem 4.2}, we obtain
\begin{align} \label{eq 4.2}
q &= -(\cot\theta)A^{\Sigma}(\nu,\nu) + \frac{1}{\sin\theta}A^{\partial M}(\bar{\nu},\bar{\nu}) \\
&= \frac{1}{\sin\theta}\left(H^{\partial M} - (\cos\theta)\Theta^{+} + (\cos\theta){\rm tr}_{\Sigma}A^{M} - (\sin\theta)H^{\partial\Sigma}\right). \nonumber
\end{align}

We now extend the free boundary result (Theorem \ref{main thm 1}) concerning positive scalar curvature with minimal boundary to the capillary boundary setting.

\begin{proof of main thm 3}
Let $\psi > 0$ be an eigenfunction of $L_{1}$ satisfying the Robin boundary condition $B_{2}\psi = 0$. Since $\Sigma$ is a $\mathfrak{g}$-stable hypersurface with capillary boundary, we have
\begin{align} \label{eq 4.3.3}
\begin{cases}
&L_{1}\psi = -\Delta_{\Sigma}\psi + \tilde{Q}\psi \geq 0 \,\, {\rm on} \,\, \Sigma, \\
&B_{2}\psi = \frac{\partial\psi}{\partial\nu} - \left(q - \langle X, \nu \rangle\right)\psi = 0 \,\, {\rm along} \,\, \partial\Sigma.
\end{cases}
\end{align}
By (\ref{eq 4.2}), we obtain
\begin{align} \label{eq 4.3}
q - \langle X, \nu \rangle = \frac{1}{\sin\theta}\left(H^{\partial M} - (\cos\theta)h + (\cos\theta){\rm tr}_{\Sigma}A^{M} - (\sin\theta)H^{\partial\Sigma} - (\sin\theta)A^{M}(N,\nu)\right).
\end{align}
We claim that 
\begin{align} \label{eq 4.4}
(\cos\theta){\rm tr}_{\Sigma}A^{M} - (\sin\theta)A^{M}(N,\nu) = (\cos\theta){\rm tr}_{\partial M}A^{M} + (\sin\theta)A^{M}(\bar{N},\bar{\nu}).
\end{align}
To verify this, let $\{e_{i}\}_{1 \leq i \leq n}$ be an orthonormal frame of $T\Sigma$ such that $e_{n} = \nu$. Then
\begin{align*}
(\cos\theta){\rm tr}_{\Sigma}A^{M} - (\sin\theta)A^{M}(N,\nu) =& (\cos\theta)\sum^{n-1}_{i=1}A^{M}(e_{i}, e_{i}) + (\cos\theta)A^{M}(\nu,\nu) \\
&- (\sin\theta)A^{M}(N,\nu) \\
=& (\cos\theta)\sum^{n-1}_{i=1}A^{M}(e_{i}, e_{i}) + A^{M}(\bar{\nu}, \nu) \\
=& (\cos\theta)\sum^{n-1}_{i=1}A^{M}(e_{i}, e_{i}) + (\cos\theta)A^{M}(\bar{\nu}, \bar{\nu}) \\
&+ (\sin\theta)A^{M}(\bar{\nu}, \bar{N}) \\
=& (\cos\theta){\rm tr}_{\partial M}A^{M} + (\sin\theta)A^{M}(\bar{\nu},\bar{N}),
\end{align*}
which proves the claim (see also Appendix A in \cite{CH4}). Substituting (\ref{eq 4.4}) into (\ref{eq 4.3}) and using the TDBEC (\ref{TDBEC}), we obtain
\begin{align} \label{eq 4.6}
q - \langle X, \nu \rangle \geq -\frac{1}{\sin\theta}\left((\cos\theta)h + (\sin\theta)H^{\partial\Sigma} \right) \geq -H^{\partial\Sigma},
\end{align}
where, in the last inequality, we used the assumptions $h \geq 0$ and $\theta \in [\frac{\pi}{2}, \pi)$. Consequently, the eigenvalue problem (\ref{eq 4.3.3}) reduces to
\begin{align*}
\begin{cases}
&-\Delta_{\Sigma}\psi + \tilde{Q}\psi \geq 0 \,\, {\rm on} \,\, \Sigma, \\
&\frac{\partial\psi}{\partial\nu} + H^{\partial\Sigma}\psi \geq 0 \,\, {\rm along} \,\, \partial\Sigma.
\end{cases}
\end{align*}
The remainder of the proof proceeds exactly as in the proof of Theorem \ref{main thm 1}, and we therefore omit the details.

\end{proof of main thm 3}

Finally, we extend the rigidity result for $\mathfrak{g}$-stable surfaces with free boundary to the case of $\mathfrak{g}$-stable surface with capillary boundary.

\begin{proof of main thm 4}
Since $\Sigma$ is a $\mathfrak{g}$-stable surface with capillary boundary, by (\ref{lem 4.1}) in Lemma \ref{thm 4.1} we have
\begin{align*}
0 \leq \int_{\Sigma} \left(\vert \nabla^{\Sigma}\psi \vert^{2} + \tilde{Q}\psi^{2}\right) dv - \int_{\partial\Sigma}\left(q - \langle X,\nu \rangle\right) ds,
\end{align*}
where $\tilde{Q} = \frac{R^{\Sigma}}{2} - \tilde{P}$ and $\tilde{P} = \left(\mu + J(N)\right) + \frac{\vert \chi \vert^{2}}{2} - h\left(\tau - \frac{h}{2}\right)$. Taking $\psi = 1$, we obtain
\begin{align*}
0 \leq \int_{\Sigma} K^{\Sigma} dv - \int_{\Sigma}\left(\mu - \vert J \vert\right) dv - \int_{\Sigma}\left(\tau - \frac{h}{2}\right)h \, dv + \int_{\partial\Sigma} k^{\partial\Sigma} ds,
\end{align*}
where we used \eqref{eq 4.6} to estimate the boundary term. By our assumptions and the Gauss-Bonnet theorem, this yields
\begin{align} \label{eq 4.7}
C_{1}A(\Sigma) \leq 2\pi\chi(\Sigma).
\end{align}
This inequality implies that $\Sigma$ is topologically a disk. Moreover, if the equality holds in (\ref{eq 4.7}), then all of the above inequalities become equalities. Consequently, properties (1)-(4) in Theorem \ref{main thm 4} follow, by arguments entirely analogous to those in the proof of Theorem \ref{main thm 2}.

\end{proof of main thm 4}

\end{document}